\def\minus{\smallsetminus}
\def\oW{\overline{W}}
\def\p{\partial}
\def\CP{\mathbb{C}P}
\def\CX{\mathbb{C}}
\def\RE{\mathbb{R}}

\def\HH{\eus{H}}
\def\HHbd{\eus{H}_{\mathrm{bd}}}
\def\LL{\eus{L}}
\def\we{\wedge}
\def\bM{\mathbf{M}}
\def\oI{\overline{I}}
\def\oJ{\overline{J}}

\def\oW{\overline{W}}
\def\oZ{\overline{Z}}
\def\tC{\tilde{C}}
\def\tS{\tilde{S}}
\def\RR{\overset{\circ}{R}}

\documentclass[reqno]{amsart}
\usepackage{amssymb,euscript}
\newcommand{\eus}{\EuScript}

\newcommand{\Diff}{\mbox{\sf Diff}}
\newcommand{\Metr}{\mbox{\sf Metr}}
\DeclareMathOperator{\Ker}{Ker}

\DeclareMathOperator{\const}{const}
\DeclareMathOperator{\dist}{dist}
\DeclareMathOperator{\Image}{Im}

\DeclareMathOperator{\Ric}{Ric}

\DeclareMathOperator{\tr}{tr}
\DeclareMathOperator{\vol}{vol}
\DeclareMathOperator{\Sym}{Sym}

\begin{document}
\renewcommand{\a}{\alpha}
\renewcommand{\b}{\beta}
\renewcommand{\d}{\delta}
\newcommand{\D}{\Delta}
\newcommand{\e}{\varepsilon}
\newcommand{\g}{\gamma}
\newcommand{\G}{\Gamma}
\newcommand{\la}{\lambda}
\newcommand{\La}{\Lambda}
\newcommand{\n}{\nabla}
\newcommand{\var}{\varphi}
\newcommand{\s}{\sigma}
\newcommand{\Sig}{\Sigma}
\renewcommand{\t}{\tau}
\renewcommand{\th}{\theta}
\renewcommand{\O}{\Omega}
\renewcommand{\o}{\omega}
\newcommand{\z}{\zeta}

\newcommand{\ben}{\begin{enumerate}}
\newcommand{\een}{\end{enumerate}}
\newcommand{\be}{\begin{equation}}
\newcommand{\ee}{\end{equation}}
\newcommand{\bea}{\begin{eqnarray}}
\newcommand{\eea}{\end{eqnarray}}
\newcommand{\bc}{\begin{center}}
\newcommand{\ec}{\end{center}}

\newcommand{\IR}{\mbox{I \hspace{-0.2cm}R}}
\newcommand{\IN}{\mbox{I \hspace{-0.2cm}N}}

\newtheorem{thm}{Theorem}[section]
\newtheorem{cor}[thm]{Corollary}
\newtheorem{lem}[thm]{Lemma}
\newtheorem{prop}[thm]{Proposition}
\newtheorem{ax}{Axiom}
\newtheorem{conj}[thm]{Conjecture}

\theoremstyle{definition}
\newtheorem{defn}{Definition}[section]

\theoremstyle{remark}
\newtheorem{rem}{\rm\bfseries{Remark}}[section]
\newtheorem*{notation}{Notation}

\title[asymptotically cylindrical Calabi--Yau manifolds]{Ricci-flat
  deformations of asymptotically cylindrical Calabi--Yau manifolds}
\author[KOVALEV]{Alexei Kovalev
}
\address{DPMMS, University of Cambridge, Centre for Mathematical Sciences,
  Wilberforce Road, Cambridge CB3 0WB, England}
\email{a.kovalev@dpmms.cam.ac.uk}
\thanks{In S.~Akbulut, T.~\"Onder, and R.J. Stern, editors, {\em Proceedings
  of G\"okova Geometry-Topology Conference}, pages 137--153.
  International Press, 2006.}

\begin{abstract}
We study a class of asymptotically cylindrical Ricci-flat K\"ahler metrics
arising on quasiprojective manifolds. Using the Calabi--Yau geometry and
analysis and the Kodaira--Kuranishi--Spencer theory and building up on results
of N.Koiso, we show that under rather general hypotheses any local
asymptotically cylindrical Ricci-flat deformations of such metrics are again
K\"ahler, possibly with respect to a perturbed complex structure. We also find
the dimension of the moduli space for these local deformations. In the
class of asymptotically cylindrical Ricci-flat metrics on $2n$-manifolds, the
holonomy reduction to $SU(n)$ is an open condition.
\end{abstract}

\maketitle

Let $M$ be a compact smooth manifold with integrable complex structure~$J$
and $g$ a Ricci-flat K\"ahler metric with respect to~$J$. A theorem due to
N.Koiso~\cite{koiso} asserts that if the deformations of the complex
structure of~$M$ are unobstructed then the Ricci-flat K\"ahler metrics
corresponding to the nearby complex structures and K\"ahler classes fill in
an open neighbourhood in the moduli space of Ricci-flat metrics on~$M$. The
proof of this result relies on Hodge theory and 
Kodaira--Spencer--Kuranishi theory and Koiso also found the dimension of
the moduli space.

The purpose of this paper is to extend the above result to a class of
complete Ricci-flat K\"ahler manifolds with asymptotically cylindrical
ends (see \S\ref{one} for precise definitions). A suitable version of
Hodge theory was developed as part of elliptic theory for asymptotically
cylindrical manifolds in \cite{LM,lockhart,MP,tapsit}. A complex manifold 
underlying an asymptotically cylindrical Ricci-flat K\"ahler manifold
admits a compactification by adding a `divisor at infinity'. There is
an extension of Kodaira--Spencer--Kuranishi theory for this class of
non-compact complex manifolds using the cohomology of logarithmic
sheaves~\cite{kawamata}.
On the other hand, manifolds with asymptotically cylindrical ends appear as
an essential step in the gluing constructions of compact manifolds endowed
with special Riemannian structures. In particular, the Ricci-flat K\"ahler
asymptotically cylindrical manifolds were prominent in~\cite{kov} in the
construction of compact 7-dimensional Ricci-flat manifolds with special
holonomy $G_2$.

We introduce the class of Ricci-flat K\"ahler asymptotically cylindrical
manifolds in~\S\ref{one}, where we also state our first main
Theorem~\ref{main} and give interpretation in terms of special holonomy.
We review basic facts about the Ricci-flat deformations in~\S\ref{two}.
\S\S 3--5 contain the proof of Theorem~\ref{main} and our second main
result Theorem~\ref{dimension} on the dimension of the moduli space for the 
Ricci-flat asymptotically cylindrical deformations of a Ricci-flat K\"ahler
asymptotically cylindrical manifold. Some examples (motivated
by~\cite{kov}) are considered in~\S\ref{examples}.

\section{Asymptotically cylindrical manifolds}
\label{one}

A non-compact Riemannian manifold $(M,g)$ is called {\em asymptotically
cylindrical} with cross-section~$Y$ if\\[3pt]
(1) $M$ can be decomposed as a union $M=M_{\text{cpt}}\cup_Y M_e$ of a
compact manifold $M_{\mathrm{cpt}}$ with boundary $Y$ and an end $M_e$
diffeomorphic to half-cylinder $[1,\infty)\times Y$, the two pieces attached
via $\p M_{\mathrm{cpt}}\cong \{1\}\times Y$, and\\[3pt]
(2) The metric $g$ on~$M$ is asymptotic, along the end, to a product
cylindrical metric $g_0=dt^2+g_Y$ on $[1,\infty)\times Y$,
$$
\lim_{t\to\infty}(g-g_0)=0,\qquad
\lim_{t\to\infty}\n_0^k g=0,\quad k=1,2,\ldots,
$$
where $t$ is the coordinate on~$[1,\infty)$ and $\n_0$ denotes the
Levi--Civita connection of~$g_0$.

Note that the cross-section $Y$ is always a compact manifold. We shall
sometimes assume that $t$ is extended to a smooth function defined on
all of~$M$, so that $t\ge 1$ on the end and $0\le t\le 1$ on the compact
piece of~$M_{\mathrm{cpt}}$.
        \begin{rem}\label{bmet}
Setting $x=e^{-t}$, one can attach to $M$ a copy of $Y$ corresponding to
$x=0$ and obtain a compactification $\bM=M\cap Y$ `with boundary at
infinity'. Then $x$ defines a normal coordinate near the boundary
of~$\bM$. The metric $g$ is defined on the interior of $\bM$ and blows up
in a particular way at the boundary,
        \begin{equation}\label{exbm}
g=\bigl(\frac{dx}{x}\bigr)^2+\tilde{g},
        \end{equation}
for some semi-positive definite symmetric bilinear $\tilde{g}$ smooth
on~$M$ and continuous on~$\bM$, such that $\tilde{g}|_{x=0}=g_Y$. Metrics
of this latter type are called `exact $b$-metrics' and are studied
in~\cite{tapsit}.
        \end{rem}
Our main result concerns a K\"ahler version of the asymptotically
cylindrical Riemannian manifolds which we now define.
Suppose that $M$ has an integrable complex
structure~$J$ and write $Z$ for the resulting complex manifold.
The basic idea is to replace a real parameter $t$ along the
cylindrical end by a complex parameter, $t+i\theta$ say, where
$\theta\in S^1$. Thus in the complex setting the asymptotic model for a
cylindrical end of $Z$ takes a slightly special form
$\RE_{>0}\times S^1\times D$, for some compact complex manifold~$D$. 
Respectively, the normal coordinate $x=e^{-t}$ becomes the real part of a
holomorphic local coordinate $z=e^{-t-i\theta}$ taking values in the
punctured unit disc $\Delta^*=\{0<|z|<1\}\subset\CX$. It follows that the
complex structure on the cylindrical end is asymptotic to the product
$\Delta^*\times D$ and the complex manifold~$Z$ is {\em compactifiable},
$Z=\oZ\minus D$, where $\oZ$ is a compact complex manifold of the same
dimension as~$Z$ and $D$ is a complex submanifold of codimension~1 in~$\oZ$
with holomorphically trivial normal bundle~$N_{D/\;\oZ}$.

The local complex coordinate $z$ on~$\oZ$ vanishes to order one precisely
on~$D$ and a tubular neighbourhood $\oZ_e=\{|z|<1\}$ is a local deformation
family for~$D$,
        \be\label{locfibr}
\pi:\oZ_e\to\Delta,\qquad
D=\pi^{-1}(0),
        \ee
where $\pi$ denotes the holomorphic map defining the coordinate~$z$.
Note that the cylindrical end $Z_e=\oZ_e\minus D$ is diffeomorphic (as a
real manifold) but not in general biholomorphic to $\RE_{>0}\times
S^1\times D$ as the complex structure on the fibre~$\pi^{-1}\{z\}$ depends
on~$z$.
        \begin{rem}
If $H^{0,1}(\oZ)=0$ then the local map~\eqref{locfibr}
extends to a holomorphic fibration $\oZ\to\CX P^1$ (cf.\
\cite[pp.34--35]{GH}).
        \end{rem}
A product K\"ahler metric, with respect to a product complex structure on
$\RE\times S^1\times D$, has K\"ahler form $a^2 dt\we d\theta+\o_D$,
where $\o_D$ is a K\"ahler form on~$D$ and $a$ is a positive
function of $t,\theta$. We shall be interested in the situation when the
product K\"ahler metric is {\em Ricci-flat}; then $a$ is a constant and can
be absorbed by rescaling the variable $t$.

We say that a K\"ahler metric on~$Z$ is {\em asymptotically cylindrical}
if its K\"ahler form $\o$ can be expressed on the end
$Z_e=\oZ_e\minus D\subset Z$ as
$$
\o|_{Z_e}=\o_D + dt\we d\theta+\e,
$$
for some closed form $\e\in\O^2(Z_e)$
decaying, with all derivatives, to zero uniformly on
$S^1\times D$ as $t\to\infty$. An asymptotically cylindrical K\"ahler
metric defines an asymptotically cylindrical Riemannian
metric on the underlying real manifold.

We shall sometimes refer to K\"ahler metrics by their K\"ahler forms.
        \begin{prop}\label{exp.decay}
Let $Z$ be a compactifiable complex manifold as defined above.
If $\o$ is an asymptotically cylindrical K\"ahler metric on~$M$ then
the decaying term on~$Z_e$ is {\em exact},
        \be\label{strong}
\o|_{Z_e}=\o_D+ dt\we d\theta+d\psi.
        \ee
        \end{prop}
        \begin{pf}
We can write $\e=\e_0(t)+dt\we\e_1(t)$, where $\e_0(t),\e_1(t)$
are 1-parameter families of, respectively, 2-forms and 1-forms on the
cross-section $S^1\times D$. As $\e$ is closed, $\e_0(t)$ must be
closed for each~$t$ and $\frac{\p}{\p t}\e_0(t)=d_{S^1\times D}\e_1(t)$.
As $\e_1$ decays exponentially fast, we have
$\e_0(t)=\int_\infty^t d_{S^1\times D}\e_1(s)ds$ and the integral
converges absolutely. So we can write
$$
\e = d_{S^1\times D} \int_\infty^t \e_1(s)ds + dt\we\e_1(t)
$$
which is an exact differential of a 1-form $\psi=-\int_t^\infty \e_1(s)ds$
on~$Z_e$.
        \end{pf}
Recall that by Yau's solution of the Calabi conjecture a compact K\"ahler
manifold admits Ricci-flat K\"ahler metrics if and only if its first Chern
class vanishes~\cite{yau}. Moreover, the Ricci-flat K\"ahler metric is
uniquely determined by the cohomology class of its K\"ahler
form. Ricci-flat K\"ahler manifolds are sometimes called Calabi--Yau
manifolds.
        \begin{rem}
There is an alternative way to define the Calabi--Yau manifolds using
the holonomy reduction. The holonomy group of a Riemannian
$2n$-manifold is the group of isometries of a tangent space
generated by parallel transport using the Levi--Civita connection over
closed paths based at a point. The holonomy group can be identified with a
subgroup of $SO(2n)$ if the manifold is orientable. If the holonomy of
a Riemannian $2n$-manifold is contained in $SU(n)\subset SO(2n)$ then
the manifold has an integrable complex structure~$J$, so that with respect
to~$J$ the metric is Ricci-flat K\"ahler. The converse is in general
not true unless the manifold is simply-connected.
        \end{rem}

A version of the Calabi conjecture for asymptotically
cylindrical K\"ahler manifolds is proved in \cite[Thm.~5.1]{TY1} and
\cite[\S\S 2--3]{kov}. It can be stated as the following.
        \begin{thm}\label{bCY}{\em (cf.~\cite[Thms. 2.4 and 2.7]{kov})}
Suppose that $Z=\oZ\minus D$ is a compactifiable complex $n$-fold as defined
above, such that $D$ is an anticanonical divisor on~$\oZ$ and the normal
bundle of $D$ is holomorphically trivial and $b^1(\oZ)=0$.
Let $\overline{g}$ be a K\"ahler metric on~$\oZ$ and denote by $g_D$ the
Ricci-flat K\"ahler metric on~$D$ in the K\"ahler class defined by the
embedding in~$\oZ$.

Then $Z=\oZ\minus D$ admits a complete Ricci-flat K\"ahler metric $g_Z$. The
K\"ahler form of $g_Z$ can be written, on the cylindrical end of~$Z$, as
in~\eqref{strong} with $\o_D$ the K\"ahler form of~$g_D$.

If, in addition, $\oZ$ and $D$ are simply-connected and there is a closed
real 2-dimensional submanifold of $\oZ$ meeting $D$ transversely with
non-zero intersection number then the holonomy of $g$ is $SU(n)$.
        \end{thm}
Note that an anticanonical divisor $D$ admits Ricci-flat K\"ahler metrics as
$c_1(D)=0$ by the adjunction formula.
The result in~\cite{kov} is stated for threefolds, but the proof
generalizes to an arbitrary dimension by a change of notation.
We consider examples arising by application of the above theorem
in~\S\ref{examples}. 
A consequence of the arguments in~\cite{kov} is that
if an asymptotically cylindrical K\"ahler metric $\o$ is Ricci-flat 
then the \mbox{1-form} $\psi\in\O^1(M_e)$ in~\eqref{strong} can
be taken to be decaying, with all derivatives, at an {\em exponential} rate
$O(e^{-\la t})$ as $t\to\infty$, for some $0<\la<1$ depending
on~$g_D$. Furthermore, if $\o$ and $\tilde\o$ are
asymptotically cylindrical Ricci-flat metrics on~$Z$ such that
$\tilde\o=\o+i\p\bar\p u$ for some $u\in C^\infty(Z)$ decaying to zero as
$t\to\infty$ then $\o=\tilde\o$ \cite[Propn.~3.11]{kov}.

Let $(M,g)$ be an asymptotically cylindrical Riemannian manifold. A local
deformation $g+h$ of $g$ is given by a field of symmetric bilinear forms
satisfying $|h|_g<1$ at each point, so that $g+h$ is a well-defined metric.
Suppose that $g+h$ is asymptotically cylindrical. Then there is a
well-defined symmetric bilinear form $h_Y$ on~$Y$ obtained as the
limit of $h$ as $t\to\infty$ and $h_Y$ is a deformation of the limit $g_Y$
of $g$, in particular $|h_Y|_{g_Y}<1$. The $h_Y$ defines via  the obvious
projection $\RE\times Y\to Y$ a $t$-independent symmetric bilinear
form on the cylinder, still denoted by~$h_Y$. Let $\rho:\RE\to[0,1]$ denote
a smooth function, such that $\rho(t)=1$, for $t\ge 2$, and
$\rho(t)=0$, for $t\le 1$. 
In view of the remarks in the previous paragraph we shall be interested in
the class of metrics which are asymptotically cylindrical at an exponential
rate and deformations $h$ satisfying $h-\rho h_Y=e^{-\mu t}\tilde{h}$ for
some $\mu>0$ and a bounded $\tilde{h}$.
Given an exponentially asymptotically cylindrical metric $g$, a
deformation $g+h$ `sufficiently close' to $g$ is understood in
the sense of sufficiently small Sobolev norms of $\tilde{h}$ and
$h_Y$, where the Sobolev norms are chosen to dominate the
uniform norms on $M$ and $Y$, respectively.

We now state our first main result in this paper.
        \begin{thm}\label{main}
Let $W=\oW\minus D$, where $\oW$ is a compact complex manifold and $D$ is
a smooth anticanonical divisor on~$\oW$ with holomorphically trivial normal
bundle. Let $g$ be an asymptotically cylindrical Ricci-flat K\"ahler metric
on~$W$. Suppose that all the compactifiable infinitesimal deformations of
the complex manifold $W$ are integrable (arise as tangent vectors to paths
of deformations).

Then any Ricci-flat asymptotically cylindrical metric on~$W$ sufficiently
close to~$g$ is K\"ahler with respect to some compactifiable deformation of
the complex structure on~$W$.
        \end{thm}
The additional conditions for the holonomy reduction given in
Theorem~\ref{bCY} are topological and we deduce from Theorem~\ref{main}.
        \begin{cor}
Assume that $W=\oW\minus D$ satisfies the hypotheses of Theorem~\ref{main}.
Suppose further that $W$, $\oW$, and $D$ are simply-connected and so the
metric $g$ has holonomy $SU(n)$, $n=\dim_\CX W$. Then any Ricci-flat
asymptotically cylindrical metric on~$W$ close to~$g$ also has
holonomy~$SU(n)$. 
        \end{cor}
Our second main result determines the
dimension of the moduli space of the asymptotically cylindrical Ricci-flat
K\"ahler metrics and is given by Theorem~\ref{dimension} below.

\section{Infinitesimal Ricci-flat deformations}
\label{two}

Before dealing with the moduli of asymptotically cylindrical Ricci-flat
metrics we recall, in summary, some results on the moduli problem for
the Ricci-flat metrics on a compact manifold. The case of a compact
manifold is standard and further details can be found
in~\cite[Ch.~12]{besse} and references therein.

A natural symmetry group of the equation $\Ric(g)=0$ for a metric $g$ on a
compact manifold~$X$ is the group $\Diff X$ of diffeomorphisms of~$X$.
It is also customary to identify a metric $g$ with $a^2g$, for
any positive constant $a$. This is equivalent to considering only the metrics
such that $X$ has volume~1. The moduli space of Ricci-flat metrics on~$X$ is
defined as the space of orbits of all the solutions $g$ of $\Ric(g)=0$ in the
action of $\Diff(X)\times\RE_{>0}$,
$$
g\mapsto a^2\phi^*g, \qquad \phi\in\Diff X,\; a>0,
$$
or, equivalently, the space of all ($\Diff X$)-orbits of the solutions of
$\Ric(g)=0$ such that $\vol_g(X)=1$.
The tangent space at~$g$ to an orbit of~$g$ under the action
of $\Diff X$ is the image of the first order linear differential operator
        \be\label{dstar}
\delta_g^*:V^\flat\in\O^1(X)\to
{\textstyle\frac12} \LL_Vg \in\Sym^2T^*X,
        \ee
where $\LL$ denotes the Lie derivative.
The operator $\delta_g^*$ may be equivalently expressed as the symmetric
component of the Levi--Civita covariant derivative
$\n_g:\O^1(X)\to\O^1\otimes\O^1(X)$, for the metric~$g$,
        \be\label{symc}
\n_g\eta=\delta_g^*\eta+\frac12 d\eta,\quad
\eta\in\Omega^1(X).
        \ee
The $L^2$ formal adjoint of $\delta_g^*$ is therefore given by
$$
\delta_g:h\in\Sym^2T^*X\to \n_g^*h\in\O^1(X).
$$
The operator $\delta_g^*$ is overdetermined-elliptic with
finite-dimensional kernel and closed image and there is an $L^2$-orthogonal
decomposition
$$
\Sym^2T^*X=\Ker\delta_g\oplus\Image\delta_g^*.
$$
The equation $\delta_gh=0$ defines a local transverse slice
for the action of $\Diff(X)$.

The infinitesimal Ricci-flat deformations $h$ of a Ricci-flat~$g$
preserving the volume are obtained by linearizing the equation $\Ric(g+h)=0$
at $h=0$, imposing an additional condition $\int_X\tr_g h\:\nu_g=0$, where
$\nu_g$ is the volume form of~$g$.
By a theorem of Berger and Ebin, the space of infinitesimal Ricci-flat
deformations of~$g$ is given by a system of linear PDEs
        \be\label{BEcomp}
(\n^*_g\n_g - 2\RR_g)h=0,\qquad
\delta_g h =0,\qquad
\tr_g h=0.
        \ee
Here $\RR_g$ is a linear map induced by the Riemann curvature
and acting on symmetric bilinear forms,
$\RR_gh(X,Y)=\sum_ih(R_g(X,e_i)Y,e_i)$ \ ($e_i$ is an orthonormal basis).
The first equation in~\eqref{BEcomp} is elliptic and so the solutions
of~\eqref{BEcomp} form a finite-dimensional space.

Suppose that every infinitesimal deformation~$h$ satisfying~\eqref{BEcomp}
arises as the tangent vector to a path of Ricci-flat metrics. Then it turns
out that a neighbourhood of~$g$ in the moduli space of Ricci-flat metrics
on~$X$ is diffeomorphic to the quotient of the solutions space
of~\eqref{BEcomp} by a finite group. This finite group depends on the
isometry group of~$g$ and the moduli space is an orbifold of dimension
equal to the dimension of the solution space of~\eqref{BEcomp}.

Now suppose that the manifold $X$ has an integrable complex
structure, $J$ say, and the Ricci-flat metric $g$ on $X$ is K\"ahler,
with respect to~$J$.
Then any deformation $h$ of~$g$ may be written as a sum $h=h_+ + h_-$
of Hermitian form $h_+$ and skew-Hermitian form $h_-$ defined by the
conditions $h_\pm(Jx,Jy)=\pm h_\pm(x,y)$. Furthermore, the operator
$\n^*_g\n_g - 2\RR_g$ preserves the subspaces of Hermitian and
skew-Hermitian forms.

The skew-Hermitian forms $h_-$ may be identified, via
        \be\label{iform}
g(x,Iy)=h_-(x,Jy),
        \ee
with the symmetric real endomorphisms $I$ satisfying $IJ+JI=0$.
Thus $J+I$ is an almost complex structure and the endomorphism $I$ may be
regarded as a $(0,1)$-form with values in the holomorphic tangent bundle
$T^{1,0}X$. Then one has
        \be\label{igauge}
\delta_g h_- = -J \circ(\bar\p^* I).
        \ee
In particular, $\delta h_-=0$ if and only if $I$ defines an class in
$H^1(X,T^{1,0}X)$, that is $I$ defines an infinitesimal
deformation of the complex manifold~$(X,J)$ (see~\cite{ks}). 
With the help of Weitzenb\"ock formula one can replace $\n_g^*\n_g-2\RR_g$ 
by the complex Laplacian for $(0,q)$-forms with values in $T^{1,0}X$
$$
((\n^*_g\n_g - 2\RR_g)h_-)(\cdot,J\cdot)=
g(\cdot,(\Delta_{\bar\p}I)\cdot).
$$
Thus $(\n^*_g\n_g - 2\RR_g)h_-=0$ precisely when
$I\in\O^{0,1}(T^{1,0}X)$ is harmonic.

Hermitian forms $h_+$ are equivalent, with the help of the complex
structure, to the real differential (1,1)-forms
        \be\label{hform}
\psi(\cdot,\cdot)=h_+(\cdot,J\cdot).
        \ee
The Weitzenb\"ock formula yields
$$
((\n^*_g\n_g - 2\RR_g)h_+)(\cdot,J\cdot) = \Delta\psi,
$$
for a Ricci-flat metric~$g$, thus $h_+$ satisfies the first equation
in~\eqref{BEcomp} if and only if $\psi\in\O^{1,1}$ is harmonic.
The other two equations in~\eqref{BEcomp} become
        \be\label{hgauge}
\tr_g h_+ = \langle\psi,\o\rangle_g,\quad\text{and}\quad
\delta_g h_+=-d^*\psi,
        \ee
where $\omega$ denotes is the K\"ahler form of~$g$.

\section{The moduli problem and a transverse slice}

We want to extend the set-up of the moduli space for Ricci-flat
metrics outlined in \S\ref{two} to the case when $(M,g)$ is an
asymptotically cylindrical Ricci-flat manifold. 
For this, we require a Banach space completion for sections of
vector bundles associated to the tangent bundle of $M$ and we
use Sobolev spaces with exponential weights. A weighted Sobolev space
$e^{-\mu t}L^p_k(M)$ is, by definition, the space of all functions
$e^{-\mu t}f$ such that $f\in L^p_k(M)$. The norm of $e^{-\mu t}f$
in $e^{-\mu t}L^p_k(M)$ is defined to be the $L^p_k$-norm of~$f$. The
definition generalizes in the usual way to vector fields, differential forms,
and, more generally, tensor fields on~$M$ with the help of the Levi--Civita
connection. Note that if $k-\dim M/p>\ell$, for some integer $\ell\ge 0$, then
there is a bounded inclusion map between Banach spaces $L^p_k(M)\to
C^{\ell}(M)$ because $(M,g)$ is complete and has bounded
curvature~\cite[\S2.7]{aubin}.

The weighted Sobolev spaces $e^{-\mu t}L^p_k(M)$ are not quite convenient
for working with bounded sections that are asymptotically $t$-independent
but not necessarily decaying to zero on the end of~$M$. We shall use
slightly larger spaces which we call, following a prototype
in~\cite{APS}, the {\em extended weighted Sobolev spaces}, denoted
$W^p_{k,\mu}(M)$.

As before, use $Y$ to denote the cross-section of the end of~$M$.
Fix once and for all a smooth cut-off function $\rho(t)$ such that
$0\le\rho(t)\le 1$, $\rho(t)=0$ for $t\le 1$, and $\rho(t)=1$ for $t\ge 2$.
Define
$$
W^p_{k,\mu}(M)=e^{\mu t}L^p_k(M)+\rho(t)L^p_k(Y)
$$
where, by abuse of notation, $L^p_k(Y)$ in the above formula is understood
as a space of \mbox{$t$-independent} functions on the cylinder $\RE\times Y$
pulled back from~$Y$. Elements in $\rho(t)L^p_k(Y)$ are well-defined as
functions supported on the end of~$M$. The norm of $f+\rho(t)f_Y$ in
$W^p_{k,\mu}(M)$ is defined as the sum of the $e^{\mu t}L^p_k(M)$-norm of
$f$ and the $L^p_k$-norm of~$f_Y$ (where $f_Y$ is interchangeably
considered as a function on~$Y$).
More generally, the extended weighted Sobolev space of $W^p_{k,\mu}$
sections of a bundle associated to $TM$ is defined in a similar manner
using parallel transport in the $t$ direction defined by the Levi--Civita
connection.

We shall need some results of the elliptic theory and Hodge theory for an
asymptotically cylindrical manifold~$(M,g)$. The Hodge Laplacian $\Delta$
on~$M$ is an instance of an {\em asymptotically translation invariant}
elliptic operator. That it, $\Delta$ can be written locally on the end of
$M$ in the form $a(t,y,\p_t,\p_y)$, where $a$ is smooth in
$(t,y)\in\RE\times Y$ and polynomial in $\p_t,\p_y$. The coefficients
$a(t,y,\p_t,\p_y)$ have a $t$-independent asymptotic model
$a_0(y,\p_t,\p_y)$ on the cylinder $\RE\times Y$, so that
$a(t,y,\p_t,\p_y)-a_0(y,\p_t,\p_y)$ decays to zero, together with all
derivatives, as $t\to\infty$.
        \begin{prop}\label{b-hodge}
Let $(M,g)$ be an oriented asymptotically cylindrical manifold with $Y$ a
cross-section of~$M$ and let $\Delta$ denote the Hodge Laplacian on~$M$.
Then there exists $\mu_1>0$ such that for $0<\mu<\mu_1$ the following holds.

(i) The Hodge Laplacian defines bounded Fredholm linear operators
$$
\D_{\pm\mu}:e^{\pm\mu t}L^p_{k+2}\Omega^r(M)\to 
e^{\pm\mu t}L^p_{k}\O^r(M)
$$
with index, respectively, $\pm(b^r(Y)+b^{r-1}(Y))$.
The image of $\D_{\pm\mu}$ is, respectively, the subspace of the forms
in $e^{\pm\mu t}L^p_{k}\O^r(M)$ which are $L^2$-orthogonal to
the kernel of $\D_{\mp\mu}$.

(ii) Any $r$-form $\eta\in\Ker\D_g\cap e^{\mu t}L^p_{k+2}\O^r(M)$
is smooth and can be written on the end $\RE_+\times Y$ of $M$ as
        \be\label{asexp}
\eta|_{\RE\times Y}
=\eta_{00}+t\eta_{10}+dt\we(\eta_{01}+t\eta_{11})+\eta',
        \ee
where $\eta_{ij}$ are harmonic forms on $Y$ of degree $r-j$ and
the $r$-form $\eta'$ is $O(e^{-\mu_1 t})$ with all derivatives.
In particular, any $L^2$ harmonic form on~$M$ is $O(e^{-\mu_1 t})$.
The harmonic form $\eta$ is closed and co-closed precisely when
$\eta_{10}=0$ and $\eta_{11}=0$, i.e. when $\eta$ is bounded.
        \end{prop}
        \begin{pf}
For (i) see \cite{LM} or \cite{tapsit}. In particular, the last claim is
just a Fredholm alternative for elliptic operators on weighted Sobolev spaces.

The clause (ii) is an application of~\cite[Theorem~6.2]{MP}.
Cf. also \cite[Propn.~5.61 and~6.14]{tapsit} proved with an assumption that
the $b$-metric corresponding to $g$ is smooth up to and on the boundary
of~$M$ at infinity. The last claim is verified by the standard integration
by parts argument.
        \end{pf}
        \begin{cor}\label{fredh}
Assume the hypotheses and notation of Proposition~\ref{b-hodge}. Suppose
also that the metric $g$ on $M$ is asymptotic to a product cylindrical
metric on $\RE_+\times Y$ at an {\em exponential} rate $O(e^{-\mu_1 t})$.
Then for $\xi\in e^{-\mu t}L^p_{k}\Omega^r(M)$, the equation $\Delta\eta=\xi$
has a solution
$\eta\in e^{-\mu t}L^p_{k+2}\Omega^r(M)+\rho(t)(\eta_{00}+dt\we\eta_{01})$ 
if and only if $\xi$ is $L^2$-orthogonal to $\HHbd^r(M)$, where $0<\mu<\mu_1$
and $\HHbd^r(M)$ denotes the space of bounded harmonic $r$-forms on~$M$.
        \end{cor}
        \begin{pf}
The hypotheses on~$g$ and $\mu$ implies that the Laplacian defines a
Fredholm operator
        \be\label{delplus}
e^{-\mu t}L^p_{k+2}\Omega^r(M)+\rho(t)(\eta_{00}+dt\we\eta_{01})
\to e^{-\mu t}L^p_{k+2}\Omega^r(M).
        \ee
It follows from Proposition~\ref{b-hodge} that the index of~\eqref{delplus}
is zero and the kernel is $\HHbd^r(M)$. Further, if
$\xi\in e^{-\mu t}L^p_{k+2}\Omega^r(M)+\rho(t)(\eta_{00}+dt\we\eta_{01})$
then we find from~\eqref{asexp} that $d\xi$ and $d^*\xi$ decay to zero as
$t\to\infty$. Recall from Proposition~\ref{b-hodge} that any bounded
harmonic form is closed and co-closed and then the standard Hodge theory
argument using integration by parts is valid and shows that the image
of~\eqref{delplus} is $L^2$-orthogonal to $\xi\in\HHbd^r(M)$.
But as the codimension of the image of~\eqref{delplus} is equal to
$\dim\HHbd^r(M)$ the image must be precisely the $L^2$-orthogonal
complement of $\HHbd^r(M)$ in $e^{-\mu t}L^p_{k+2}\Omega^r(M)$.
        \end{pf}
For an asymptotically cylindrical $n$-dimensional manifold $(M,g)$,
let $\Diff_{p,k,\mu}M$ (where $k-n/p>1$, $0<\mu<\mu_1$) denote the
group of locally $L^p_k$ diffeomorphisms of $M$ 
generated by $\exp_V$, for all vector fields $V$ on~$M$ that can be written
as $V=V_0+\rho(t)V_Y$, where $V_0\in e^{-\mu t}L^p_k$ and a
\mbox{$t$-independent} $V_Y$ is defined by a Killing field for $g_Y$.
Respectively, $V_Y^\flat$ is defined by a harmonic 1-form on~$Y$
(cf.~\eqref{rel} below). Also
require that $V$ has a sufficiently small $C^1$ norm on~$M$, so that that
$\exp_V$ is a well-defined diffeomorphism.
Denote by $\Metr_{p,k,\mu}(g)$ (where $k-n/p>0$, $0<\mu<\mu_1$) the space of 
deformations $h+\rho(t)h_Y$ of~$g$ where $h\in e^{-\mu t}L^p_k$,
$|h|_g<1$ at
every point of~$M$, and $\delta_{g_Y} h_Y=0$, $\tr_{g_Y} h_Y=0$.
($g_Y$ is the limit of~$g$ as defined in~\S\ref{one}.)
Then $\Diff_{p,k+1,\mu}M$ acts on $\Metr_{p,k,\mu}(g)$ by pull-backs
and the linearization of the action is given by the operator $\delta_g^*$
on weighted Sobolev spaces,
        \be
\d_g^*: e^{-\mu t}L^p_{k+1}\O^1(M)+\rho(t)\HH^1(Y)\to
e^{-\mu t}L^p_k\Sym^2T^*M.
        \ee
It will be convenient to replace the last two equations
in~\eqref{BEcomp} and instead use another local slice
equation for the action of $\Diff_{p,k,\mu}M$
$$
\d_g h+\frac12 d\tr_g h=0.
$$
A transverse slice defined by the operator $\d_g+\frac12 d\tr_g$ was
previously used for different classes of complete non-compact manifolds
in~\cite[I.1.C and I.4.B]{biq-ast}.
The operator $\d_g+\frac12 d\tr_g$ satisfies a useful relation:
        \be\label{rel}
(2\d_g+d\tr_g)\d_g^*=2\n^*_g\n_g-\n^*_g d + d\tr_g\d_g^*=
2\n_g^*\n_g-d^*_g d-dd^*_g=\D_g,
        \ee
where $\Delta_g$ is the Hodge Laplacian and we used the Weitzenb\"ock
formula for 1-forms on a Ricci-flat manifold in the last equality.
        \begin{prop}\label{fredh-alt}
Assume that $Y$ is connected and that $k-\dim M/p>1$, $0<\mu<\mu_1$,
where $\mu_1$ is defined in Propn.~\ref{b-hodge} for the Laplacian on
differential forms on~$M$.
Then there is a direct sum decomposition into closed subspaces
        \be\label{fredh.alt}
\Metr_{p,k,\mu}(g)=
\delta^*_g(e^{-\mu t}L^p_{k+1}\Omega^1(M)+\rho(t)\HH^1(Y))
\oplus \bigl(\Ker(\d_g+\frac12 d\tr_g)\cap\Metr_{p,k,\mu}(g)\bigr).
        \ee
        \end{prop}
        \begin{pf}
Any bounded harmonic 1-form on~$M$ is in
$e^{-\mu t}L^p_{k+1}\Omega^1(M)+\rho(t)\HH^1(Y)$ by\linebreak
\cite[Propn.\ 6.16 and~6.18]{tapsit} (see also Propn.~\ref{bH}
below) and because $Y$ is connected. 
For any $\eta\in e^{-\mu t}L^p_{k+1}\O^1(M)+\rho(t)\HH^1(Y)$,
$\n_g\eta$ decays on the
end of~$M$, so the standard integration by parts applies to show that
the bounded harmonic 1-forms on $M$ are parallel with respect to~$g$.
As the bounded harmonic 1-forms on~$M$ are closed we obtain using~\eqref{symc}
and~\eqref{rel} that these are in the kernel of~$\d^*$. It follows that the two
subspaces in~\eqref{fredh.alt} have trivial intersection.

By the definition of $\Metr_{p,k,\mu}(g)$ the `constant term' $h_Y$ of
$h$ satisfies $\delta_Y h_Y=0$ and $\tr_g h_Y=0$.
If $\eta\in\HHbd^1(M)$ and 
$h\in\Metr_{p,k,\mu}(g)$ then the 1-form $\langle\eta,h\rangle_g$ decays
along the end of~$M$ and we can integrate by parts
$$
\langle\eta,\d_g h+\frac12d\tr_g h\rangle_{L^2}=
\langle\d^*_g\eta,h\rangle_{L^2}+
\frac12\langle d^*\eta,\tr_g h\rangle_{L^2}=0.
$$
Thus the image $(\delta_g+\frac12d\tr_g)\Metr_{p,k,\mu}(g)$ is
$L^2$-orthogonal to $\HHbd^1$.
By Corollary~\ref{fredh} the equation
$\Delta\eta=(\delta+\frac12d\tr_g)h$ has a solution $\eta$ in
$e^{-\mu t}L^p_{k+1,\mu}\Omega^1(M)+\rho(t)\HH^1_Y$ and so
$$
\d^*_g\eta - h \in\Ker(\d_g+\frac12 d\tr_g)\cap\Metr_{p,k,\mu}(g)
$$
which gives the required decomposition $h=\d^*_g\eta+(\d^*_g\eta-h)$.
        \end{pf}
        \begin{prop}\label{gfix}
Assume that $p,k,\mu$ are as in Proposition~\ref{fredh-alt}.
Let $\tilde{g}$ an asymptotically cylindrical deformation of~$g$.
If $\tilde{h}=\tilde{g}-g\in\Metr_{p,k,\mu}(g)$ is sufficiently small in
$W^p_{k,\mu}\Sym^2T^*M$ then there exists $\phi\in\Diff_{p,k+1,\mu} M$
such that $\phi^*\tilde g=g+h$, for some $h\in\Metr_{p,k,\mu}(g)$ with
$(\delta_g+\frac12 d\tr_g)h=0$.
        \end{prop}
        \begin{pf}
If the desired $\phi$ is close to the identity then $\phi=\exp_V$ for a
vector field $V$ on~$M$ with small $e^{-\mu t}L^p_k$ norm. 
We want to show that the map
$$
T_{\mathrm{id}}(\Diff_{p,k+1,\mu}) \times 
\{h\in\Metr_{p,k,\mu}(g):(\d_g+\frac12 d\tr_g)h=0\}
\to \Metr_{p,k,\mu}(g)
$$
defined by
$$
(V,h)\mapsto \exp_V^*(g+h)-g
$$
is a onto a neighbourhood of $(0,0)$. The linearization of $(DF)_{(0,0)}$
is given by $(V,h)\mapsto\delta^*_g(V^\flat)+h$ and is surjective
by~\eqref{fredh.alt}.  By the implicit function theorem for Banach spaces,
a solution $(V,h)$ of $F(V,h)=\tilde g$ exists, whenever $\tilde g - g$ is
sufficiently small. 
        \end{pf}
Finally, we obtain the system of linear PDEs describing the infinitesimal
Ricci-flat deformations of an asymptotically cylindrical metric transverse
to the action of the diffeomorphism group on the asymptotically cylindrical
metrics.
        \begin{thm}
Suppose that $(M,g)$ is a Ricci-flat asymptotically cylindrical Riemannian
manifold, but not a cylinder $\RE\times Y$, and $g(s)$, $|s|<\e$ ($\e>0$) is
a smooth path of asymptotically cylindrical Ricci-flat metrics on~$M$ with
$g(0)=g$. Suppose also that 
$g(s)-g\in\Metr_{p,k,\mu}(g)$, with $p,k,\mu$ as in
Proposition~\ref{fredh-alt}. Then there is a smooth path
$\psi(s)\in\Diff_{p,k,\mu}M$, so that
$h=\frac{d}{ds}|_{s=0}\bigl[\psi(s)^*g(s)\bigr]$ 
satisfies the equations
        \begin{subequations}\label{inftsem}\begin{gather}
(\n^*_g\n_g - 2\RR_g)h=0,\\
(\delta_g+\frac12 d\tr_g)h=0,  \label{slice}
        \end{gather}\end{subequations}
Furthermore, if every bounded solution $h$ of~\eqref{inftsem} is the
tangent vector at~$g$ to a path of Ricci-flat asymptotically cylindrical
metrics on~$M$ then the moduli space is an orbifold. The dimension of this
orbifold is equal to the dimension of the space of the bounded on~$M$
solutions of~\eqref{inftsem}.
        \end{thm}
        \begin{pf}
Applying Proposition~\ref{gfix} for each $g(s)$, we can find a path of
diffeomorphisms in $\psi(s)\in\Diff_{p,k,\mu}M$ so that the slice
equation~\eqref{slice} holds for~$h$.

The linearization of $\Ric(g+h)=0$ in~$h$ is
$\nabla_g^*\nabla_gh-2\d^*_g\d_gh-\nabla_gd\tr_gh-2\RR_gh=0$.
which becomes equivalent to $(\n^*_g\n_g - 2\RR_g)h=0$ in view of
of~\eqref{slice} and~\eqref{symc}.

The last claim follows similarly to the case of a compact base
manifold, cf.~\cite[12.C]{besse}. It can be shown using
Proposition~\ref{fredh-alt} that the infinitesimal action of the identity
component of the group $I(M,g)$ of isometries of $g$ in $\Diff^p_{k,\mu}M$
is trivial on the slice $(\delta_g+\frac12 d\tr_g)h=0$. As $M$ is not a
cylinder, 
it has only one end~\cite{salur} and we show in Lemma~\ref{iso} below 
that $I(M,g)$ is compact. It follows that a neighbourhood of the orbit
of~$g$ in the orbit space $\Metr_{p,k,\mu}(g)/\Diff_{p,k,\mu}M$ is
homeomorphic to a finite quotient of the kernel of\linebreak
\mbox{$\delta_g+\frac12 d\tr_g$.}
        \end{pf}
        \begin{lem}\label{iso}
Let $(M,g)$ be an asymptotically cylindrical manifold with a connected
cross-section~$Y$ (that is, $M$ has only one end). Then the group $I(M,g)$
of isometries of $M$ is compact. 
        \end{lem}
        \begin{pf}
It is a well-known result the isometry group $I(M,g)$ of any
Riemannian manifold $(M,g)$ is a finite-dimensional Lie group and
if a sequence $T_k\in I(M,g)$ is such that, for some $P\in M$,
$T_k(P)$ is convergent then $T_k$ has a convergent
subsequence~\cite{myers-steenrod}.

For an  asymptotically cylindrical $(M,g)$, it is not difficult to
check that there is a choice of point $P_0$ on the end of~$M$ and
$r>0$, so that $M_{0,r}=\{P\in M:\dist(P_0,P)>r$ is connected but for any
$P_1$ such that $\dist(P_0,P_1)>3r$ the set
$M_{1,r}=\{P\in M:\dist(P_1,P)>r$ is not connected. It follows that
for any sequence $\tilde T_k\in I(M,g)$ we must have
$\dist(P_0,\tilde T_k(P_0))\le 3r$ and hence $\tilde T_k$ has a convergent
subsequence. 
        \end{pf}

\section{Infinitesimal Ricci-flat deformations of asymptotically cylindrical
K\"ahler manifolds} 

We now specialize to the {\em  K\"ahler} Ricci-flat metrics. It is
known \cite{koiso} that if an infinitesimal deformation $h$ of a Ricci-flat
K\"ahler metric on a {\em compact} manifold satisfies the Berger--Ebin
equations~\eqref{BEcomp} then the Hermitian and skew-Hermitian
components $h_+$ and $h_-$ of $h$ also satisfy~\eqref{BEcomp}. In
this section we prove a version of this result for the asymptotically
cylindrical manifolds.

        \begin{prop}\label{skew}
Let $W$ be a compactifiable complex manifold with $g$ an
asymptotically cylindrical Ricci-flat K\"ahler metric on~$W$, as
defined in~\S\ref{one}.
Suppose that an asymptotically cylindrical deformation $h\in\Metr$ of~$g$
satisfies~\eqref{inftsem}. Then the skew-Hermitian component $h_-$ of
$h$ also satisfies~\eqref{inftsem}.
        \end{prop}
        \begin{pf}
The proof uses the same ideas as in the case of for a compact
manifold (\cite[\S 7]{koiso} or \cite[Lemma~12.94]{besse}). The operator
$\n^*_g\n_g-2\RR_g$, for a 
K\"ahler metric $g$, preserves the subspaces of Hermitian and
skew-Hermitian forms, so $(\n^*_g\n_g-2\RR_g)h_-=0$. Recall from
\S\ref{two} that the latter equation implies that the form
$I\in\O^{0,1}(T^{1,0})$ corresponding to $h_-$ via~\eqref{iform} is
harmonic, $\Delta_{\bar\p}I=0$. An argument similar to that of
Proposition~\ref{b-hodge} shows that a bounded 
harmonic section~$I$ satisfies $\bar\p I=0$ which implies $\d_gh_-=0$
by~\eqref{igauge} and, further, $\d_g-\frac12 d\tr_g h_-=0$ as a
skew-Hermitian deformation $h_-$ is automatically trace-free.
        \end{pf}
        \begin{prop}
Any infinitesimal Ricci-flat asymptotically cylindrical deformation
$h\in\Metr_{p,k,\mu}(g)$ of a Ricci-flat K\"ahler asymptotically
cylindrical metric~$g$ on~$W$ is the sum of a Hermitian and a
skew-Hermitian infinitesimal deformation.

The space of skew-Hermitian infinitesimal Ricci-flat asymptotically
cylindrical deformations of~$g$ is isomorphic to the space of bounded
harmonic $(0,1)$-forms on $W$ with values in $T^{1,0}(W)$.

The space of Hermitian infinitesimal Ricci-flat asymptotically
cylindrical deformations of~$g$ is isomorphic to the orthogonal
complement of the K\"ahler form of~$g$ in the space of bounded
harmonic real $(1,1)$-forms on $W$.
        \end{prop}
        \begin{pf}
Only the last statement requires justification. Let $\o$ denote the
K\"ahler form of~$g$.

Recall from \S\ref{two} that the equation
$(\n^*_g\n_g-2\RR_g)h_+=0$ satisfied by a Hermitian
infinitesimal Ricci-flat asymptotically cylindrical deformations $h_+$
is equivalent to the condition that $\psi\in\O^{1,1}(W)$ defined
in~\eqref{hform} is harmonic $\D\psi=0$. Hence $\d_gh_+=0$
by~\eqref{hgauge} and 
Proposition~\ref{b-hodge} and so the second equation in~\eqref{inftsem}
tells us that $\langle\psi,\o\rangle_g=\const$, in view of~\eqref{hgauge}.
Considering the limit as $t\to\infty$ and the definition of
$\Metr_{p,k,\mu}(g)$ we find that the latter constant must be zero.
        \end{pf}
Thus in order to find the dimension of the space of infinitesimal
Ricci-flat deformations of an asymptotically cylindrical K\"ahler metric,
we may consider the Hermitian and a skew-Hermitian cases separately.
This is done in the next subsection.

\subsection{Bounded harmonic forms and logarithmic sheaves}

It is well-known that harmonic forms on a compact manifold are identified
with the de Rham cohomology classes via Hodge theorem. On a non-compact
manifold~$W$ one can consider the usual de Rham cohomology $H^*(W)$ and
also the de Rham cohomology $H_c^*(W)$ {\em with compact support}. The
latter is the cohomology of the de Rham complex of compactly
supported differential forms. We shall write $b^r(W)=\dim H^r(W)$ and
$b_c^r(W)=\dim H_c^r(W)$, for the respective Betti numbers. There is a
natural inclusion homomorphism $H_c^r(W)\to H^r(W)$ whose image is the
subspace of the de Rham cohomology classes representable by closed forms
with compact support; the dimension of this subspace will be denoted
by $b^r_0(W)$.
        \begin{prop}\label{bH}
Let $(W,g)$ be an oriented asymptotically cylindrical manifold. Then the space
$\HH_{L^2}(W)$ of $L^2$ harmonic $r$-forms on~$W$ has dimension
$b^r_0(W)$. The space $\HHbd(W)$ of bounded harmonic $r$-forms on~$W$ has
dimension $b^r(W) + b_c^r(W) - b^r_0(W)$.
        \end{prop}
        \begin{pf}
For the claim on $L^2$ harmonic forms see~\cite[Propn.~4.9]{APS}
or~\cite[\S 7]{lockhart}. In the case when an asymptotically cylindrical
metric~$g$ corresponds to an exact $b$-metric smooth up to the boundary
at infinity (see Remark~\ref{bmet}), the dimension of bounded harmonic forms is
a direct consequence of~\cite[Propn.~6.18]{tapsit} identifying a
Hodge-theoretic version of the long exact sequence
        \be\label{dR}
\ldots\to H^{r-1}(Y)\to H^r_c(W)\to H^r(W)\to h^r(Y)\to\ldots.
        \ee
The argument of~\cite[Propn.~6.18]{tapsit} can be adapted for
arbitrary asymptotically cylindrical metrics; the details will appear
in~\cite{mshgt}.
        \end{pf}
If $W$ is an asymptotically cylindrical K\"ahler manifold then
there is a well-defined subspace ${\HHbd^{1,1}}_{,\RE}(W)\subset
\HHbd^2(W)$ of bounded harmonic real forms of type $(1,1)$.
The bounded harmonic 2-forms in the orthogonal complement of
${\HHbd^{1,1}}_{,\RE}(W)$ are the real and imaginary parts of bounded harmonic
$(0,2)$-forms. We shall denote the complex vector space of bounded
harmonic $(0,2)$-forms on~$W$ by $\HHbd^{0,2}(W)$.

The space of bounded harmonic real (1,1)-forms on~$W$ orthogonal to the
K\"ahler form~$\o$ therefore has dimension
$b^r(W) + b_c^r(W) - b^r_0(W)-1-2\dim_\CX\HHbd^{0,2}(W)$.
\vskip 5pt

Now for the skew-Hermitian infinitesimal deformations.
Recall from \S\ref{one} that the definition of an asymptotically
cylindrical Ricci-flat K\"ahler manifold $(M,J,\omega)$ includes the
condition that a complex manifold~$W=(M,J)$ is compactifiable.
That is, there exist a compact complex $n$-fold $\oW$ and a compact
complex $(n-1)$-dimensional submanifold $D$ in $\oW$, so that $W$
is isomorphic to $\oW\minus D$.
We saw in Proposition~\ref{skew} that any skew-Hermitian Ricci-flat
asymptotically translation-invariant deformation of $\omega$ can be
expressed as a $\bar\p$- and $\bar\p^*$-closed symmetric $(0,1)$-form $I$
with values in the holomorphic tangent bundle of~$W$. A $\bar\p$- and
$\bar\p^*$-closed such $I$, not necessarily symmetric, defines an
infinitesimal deformation $J+I$ of the integrable complex
structure $J$ on~$W$.
The deformations given by skew-symmetric such forms~$I$
correspond to the bounded harmonic $(2,0)$-forms on~$W$.

Let $z$ denote a complex coordinate on $\oW$ so that $D$ is defined by the
equation $z=0$, as in~\S\ref{one}. Let $T_{\oW}$ denote the sheaf of
holomorphic local vector fields on~$\oW$. The subsheaf of the holomorphic
local vector fields whose restrictions to~$D$ are tangent to~$D$ is denoted
by $T_{\oW}(\log D)$ and called the {\em logarithmic tangent sheaf}.
The form $I$ in general has a simple pole precisely along~$D$ and defines a
class in the \v{C}ech cohomology $H^1(T_{\oW}(\log D))$.
The classical Kodaira--Spencer--Kuranishi theory of deformations of the
holomorphic structures on compact manifolds~\cite{ks} has an extension for
the compactifiable complex manifolds; the details can be found
in~\cite{kawamata}. In this latter theory, the cohomology groups
$H^i(T_{\oW}(\log D))$ have the same role as the cohomology of tangent
sheaves for the compact manifolds. In particular, the isomorphisms classes
of infinitesimal deformations of $W$ are canonically parameterized by
classes in $H^1(T_{\oW}(\log D))$. These classes arise from the actual
deformations of~$W$ is the obstruction space $H^2(T_{\oW}(\log D))$ vanishes.

Thus the space of the skew-Hermitian Ricci-flat asymptotically cylindrical
deformations $I$ of the Ricci-flat K\"ahler asymptotically cylindrical
metric $g$ on~$W$ is identified as a subspace of the infinitesimal
compactifiable deformations of~$W$. The real dimension of this subspace is
$2(\dim_\CX H^1(T_{\oW}(\log D))-\dim_\CX\HHbd^{0,2}(W))$.

\section{The asymptotically cylindrical Ricci-flat deformations}

In this section, we show that every infinitesimal Ricci-flat
deformation of an asymptotically cylindrical Ricci-flat K\"ahler
manifold is tangent to a genuine deformation. 
        \begin{thm}\label{dimension}
Let $(W,g)$ be as in Theorem~\ref{main}. Then every bounded solution $h$
of~\eqref{inftsem} arises as $h=\frac{d}{ds}|_{s=0}g(s)$ for some path of
asymptotically cylindrical Ricci-flat metrics on~$W$ with $g(0)=g$.
The moduli space of asymptotically cylindrical Ricci-flat deformations
of~$g$ is an orbifold of real dimension
$$
2\dim_\CX H^1(T_W(\log D))+b^2(W)+b^2_c(W)-b^2_0(W)-1-4\dim_\CX\HHbd^{2,0}(W).
$$
        \end{thm}
        \begin{pf}
By the hypotheses of Theorem~\ref{main}, there is a manifold $\eus{M}$ of 
small compactifiable deformations of~$W$, so that $H^1(T_{\oW}(\log D))$ is
the tangent space to~$\eus{M}$ at~$W$. The data of the compactifiable
deformations of~$W$ includes the deformations of $\oW$ \cite{kawamata}.
Let $\omega'$ be a K\"ahler metric on~$\oW$.
By the results of Kodaira and Spencer~\cite{ks}, for a family of
sufficiently small deformations~$\oI$ of a compact complex manifold~$\oW$,
there is a family of forms $\o'(\oJ+\oI)$ on~$\oW$ 
depending smoothly on~$\oI$ and such that $\o'(\oJ)=\o'$ and
$\o'(\oJ+\oI)$ defines a K\"ahler metric with respect to a perturbed
complex structure $\oJ+\oI$. Using the methods of~\cite[\S 3]{kov}, we
can construct from $\o'(\oJ+\oI)$ a smooth family $\o(J+I)$ of
asymptotically cylindrical K\"ahler metrics (not necessarily
Ricci-flat) on the respective deformations of $W=\oW\minus D$.

Consider a vector bundle $\eus{V}$ over $\eus{M}$ whose fibre over
$\oI\in\eus{M}_{\oW}$ is the space of bounded harmonic (1,1)-forms with
respect to the K\"ahler metric $\o(\oJ+\oI)$.
The task of integrating an infinitesimal Ricci-flat deformation of the
given asymptotically cylindrical K\"ahler metric $\o$ on~$W$ is expressed by
the complex Monge--Amp\'ere equation (with parameters)
for a function $u$ on~$W$
        \be\label{cMA}
(\o(J+I)+\b+i\p\bar\p u)^n-e^{f_{I,\beta}}(\o(J+I)+\b)^n=0,
        \ee
where $n=\dim_\CX W$ and $\b\in\eus{V}$ is a bounded harmonic real
(1,1)-form with respect to the K\"ahler metric $\o(J+I)$ and orthogonal to
$\o(J+I)$. The operators $\p$,$\bar\p$ in~\eqref{cMA} are those defined
by~$J+I$.

If $I=0$ and $\b=0$ then $u=0$ is a solution of~\eqref{cMA} as the metric
$\o$ is Ricci-flat. Consider the right-hand side of~\eqref{cMA} as a
function $f(I,\b,u)$ where the domain of $u$ is a version of extended
weighted Sobolev space
$E^p_{k,\mu}(W)=e^{-\mu t}L^p_{k+2}(W)+\{\rho(t)(at+b)\;|\;a,b\in\RE\}$
for a sufficiently small $\mu>0$ ($Y=S^1\times D$ is the cross-section
of~$W$ in the present case and $D$ is connected). The linearization of~$f$
in~$u$ at $u=0$ is the Laplacian for functions on the asymptotically
cylindrical K\"ahler manifold $(W,\o)$. A dimension counting argument
similar to that in Corollary~\ref{fredh} shows that this latter Laplacian
defines a {\em surjective} linear map
$E^p_{k,\mu}(W) \to e^{-\mu t}L^p_{k}(W)$.
The Laplacian has a one-dimensional kernel given by the constant
functions on~$W$, so we reduce the domain for $u$ by taking the $L^2$
orthogonal complement of the constants.
Then the implicit function theorem applies
to $f(I,\b,u)$ and defines a smooth family $u=u(I,\b)$ so that
$f(I,\b,u(I,\b))=0$ for every small $I,\b$ in the respective spaces of
bounded harmonic forms.
This defines a smooth family of Ricci-flat metrics
$\o(J+I)+\b+i\p\bar\p u(I,\b)$ tangent to the
infinitesimal deformations identified in the previous section.
        \end{pf}

\section{Examples}
\label{examples}

In this section, we consider some examples of asymptotically
cylindrical Ricci-flat K\"ahler manifolds arising by application of
Theorem~\ref{bCY} and compute the dimension of the moduli space for their
asymptotically cylindrical Ricci-flat deformations.
This is done by considering appropriate long exact sequences and
applying vanishing theorems to determine the dimensions of cohomology
groups appearing in Theorem~\ref{dimension}.

\subsection{Rational elliptic surfaces}
\label{surf}

An elliptic curve $C=\CX/\Lambda$ embeds in the complex projective plane
as a cubic curve in the anticanonical class. Choosing another non-singular
elliptic curve $C'$ in $\CP^2$ we obtain a pencil $a C+ b C'$, \
$a\! :\! b\in\CP^1$. Assuming that $C'$ is chosen generically and blowing
up the 9 intersection points $C\cap C'$ we obtain an algebraic surface
$\tilde{S}$ so that the proper transform $\tilde{C}$ of~$C$ is in the
anticanonical class, $\tilde{C}\in|-K_{\tilde{S}}|$, and $\tilde{C}$ has a
holomorphically trivial normal bundle, in particular
$\tilde{C}\cdot\tilde{C}=0$. Then, by
Theorem~\ref{bCY}, the quasiprojective surface $S=\tilde{S}\minus\tilde{C}$
has a complete Ricci-flat K\"ahler metric asymptotic to the flat metric on
the half-cylinder $\RE_{>0}\times S^1\times\CX/\Lambda$ with cross-section
a 3-dimensional torus.
Although in this example the divisor at infinity is not
simply-connected it can be easily checked that $S$ is simply-connected
and the asymptotically cylindrical Ricci-flat K\"ahler metric on~$S$
has holonomy $SU(2)$ (cf.~\cite[Theorem~2.7]{kov}). It is well-known that a
Ricci-flat K\"ahler metric on a complex surface is hyper-K\"ahler.

Furthermore, $S$ is topologically a `half of the K3 surface'
in the sense that there is an embedding of a 3-torus $T^3$ in the K3
surface so that the complement of this $T^3$ consists of two
components, each homeomorphic to $S$. From the arising Mayer--Vietoris
exact sequence, we find that $b^2(S)=b^2_c(S)=11$ using also
the Poincar\'e duality. The long exact sequence~\eqref{dR} with $W=S$
and $Y=T^3$ yields $b^2_0=8$.

As $S$ is simply-connected with holonomy $SU(2)$ there is a
nowhere-vanishing parallel (hence holomorphic) $(2,0)$-form~$\O$
on~$S$. Any other $(2,0)$-form on~$S$ can be written as $f\O$ for some
complex function ~$f$ and $f\O$ will be a bounded harmonic form if and
only if the real and imaginary parts of $f$ are bounded harmonic
functions, hence constants by the maximum principle. 
Thus $\dim_\CX \HHbd^{2,0}(S)=1$.

The dimensions of $H^1(T_S(\log\tC))$ and
$H^2(T_S(\log\tC))$ are obtained by taking the cohomology of the exact
sequences
$$
0\to T_{\tS}(-C)\to T_{\tS}(\log\tC)\to T_{\tC}\to 0
$$
and 
$$
0\to T_{\tS}(\log\tC)\to T_{\tS}\to N_{\tC/\tS}\to 0
$$
(see~\cite{kawamata}). Using Serre duality 
\cite{GH} we find that $H^2(T_{\tS}(-C))=H^0(\O^1_{\tS})=H^{0,1}(S)=0$,
hence $H^2(T_S(\log\tC))$ vanishes and the compactifiable deformations of
$S$ are unobstructed. Note that any small deformation of
$\tS$ is the blow-up of a small deformation of the cubic $\tC$ in $\CP^2$
(\cite{FN} or~\cite[Theorem~9.1]{H}).
Therefore, $\dim_\CX H^1_{\tS}=10$ and
we deduce that $\dim_\CX H^1(T_S(\log\tC))=10$.

Now by Theorem~\ref{dimension} the moduli space of
asymptotically cylindrical Ricci-flat deformations of~$S$ has
dimension 29. All these deformations are hyper-K\"ahler with
holonomy $SU(2)$.

\subsection{Blow-ups of Fano threefolds}

A family of examples of asymptotically cylindrical Ricci-flat 
K\"ahler threefolds is constructed in~\cite[\S 6]{kov} using Fano threefolds.
A Fano threefold is a non-singular complex threefold $V$ with
$c_1(V)>0$. Any Fano threefold is necessarily projective and
simply-connected. A generically chosen anticanonical divisor $D_0$ in $V$ is
a K3 surface~\cite{sh1}. Let $D_1\in|-K_V|$ be another K3 surface such that
$D_0\cap D_1=C$ is a smooth curve.

The blow-up of $V$ along~$C$ is a K\"ahler complex threefold~$(\oW,\o')$
and the proper transform~$D\subset\oW$ of $D_0$ is an anticanonical divisor
on~$\oW$ with the normal bundle of~$D$ holomorphically trivial. The
complement $W=\oW\minus D$ is simply-connected.

Thus $W$ is topologically a manifold with a cylindrical end $\RE_{>0}\times
S^1\times D$. By Theorem~\ref{bCY} $W$ admits a complete Ricci-flat
K\"ahler metric $\o$, with holonomy $SU(3)$. The metric $\o$ is asymptotic
on the end of~$W$ to the product of the standard flat metric on
$\RE_{>0}\times S^1$ and a Yau's hyper-K\"ahler metric on $D$.

By the Weitzenb\"ock formula, the Hodge Laplacian $\D$ for the $(2,0)$-forms on
a Ricci-flat K\"ahler manifold can be expressed as $\D=\n_g^*\n_g$.
The quantity $\langle\n\eta,\eta\rangle_g$ for a bounded harmonic
form~$\eta$ decays on the end of~$W$, so we can integrate by parts to show
that a bounded harmonic $(2,0)$-form is parallel. But the holonomy of the
metric $\o$ is $SU(3)$ which has no invariant elements in $\La^{2,0}\CX^3$.
Therefore, $W$ admits no parallel $(2,0)$-forms and thus no
bounded harmonic $(2,0)$-forms. 
 
The dimension of the moduli space for asymptotically cylindrical
Ricci-flat deformations of~$\o$ then becomes
$2\dim_\CX H^1(T_W(\log D))+b^2(W)+b^2_c(W)-b^2_0(W)-1$.

The dimensions of $H^i(T_W(\log D))$, $i=1,2$, are obtained from the two long
\mbox{exact} sequences similar to \S\ref{surf}. To verify that the
compactifiable deformations of $W$ are unobstructed note that
$H^2(T_{\oW})=H^1(\O^1_{\oW}(-D))=0$ by the Kodaira vanishing \mbox{theorem}
and $H^1(N_{D/\oW})=H^{1,0}(D)=0$. It is shown in \cite[\S 8]{kov} that
$b^2(W)=\rho(V)$ and $h^{2,1}(\oW)=h^{2,1}(V)+g(V)$, where
$g(V)=-K_V^3/2+1$ is the genus of $V$ and $\rho(V)$ is the Picard number.
Taking the cohomology of 
$0\to T_{\oW}(-D)\to T_{\oW}(\log D)\to T_D\to 0$
we obtain
$\dim_\CX H^1(T_W(\log D))=20+h^{2,1}(V)+g(V)-\rho(V)$.
From the long exact sequence~\eqref{dR} we find that
$b^2(W)+b^2_c(W)-b^2_0(W)=b^2(W)+1$.

Thus the dimension of the moduli space for $W$ in this example is given by 
$$
b^3(V)+2g(V)-\rho(V)+40
$$
in terms of standard invariants of the Fano threefold.
\bigskip

{\bf Acknowledgements:} The work on this paper began while the author was
visiting Universit\'e Paris~XII. I am grateful to Frank Pacard for the
invitation and many helpful discussions. I also thank Akira Fujiki and Nick
Shepherd-Barron for discussions on the logarithmic sheaves.
I would like to thank the organizers for the invitation to speak at
the G\"okova Geometry-Topology Conference, supported by T\"{U}B\.{I}TAK and
NSF.
\enlargethispage{2.ex}

\end{document}